\RequirePackage{ifpdf}
\ifpdf 
\documentclass[pdftex]{sigma}
\else
\documentclass{sigma}
\fi

\numberwithin{equation}{section}

\begin{document}

\allowdisplaybreaks

\renewcommand{\thefootnote}{$\star$}

\renewcommand{\PaperNumber}{010}

\FirstPageHeading

\ShortArticleName{$q$-Analog of Gelfand--Graev Basis for the Noncompact Quantum Algebra
$U_q(u(n,1))$}

\ArticleName{$\boldsymbol{q}$-Analog of Gelfand--Graev Basis\\ for the Noncompact Quantum Algebra
$\boldsymbol{U_q(u(n,1))}$\footnote{This paper is a contribution to the Proceedings of the XVIIIth International Colloquium on Integrable Systems and Quantum Symmetries (June 18--20, 2009, Prague, Czech Republic).  The full collection is
available at
\href{http://www.emis.de/journals/SIGMA/ISQS2009.html}{http://www.emis.de/journals/SIGMA/ISQS2009.html}}}

\Author{Raisa M. ASHEROVA~$^{\dag}$, \v{C}estm\'{\i}r BURD\'{I}K~$^{\ddag}$, Miloslav
HAVL\'{I}\v{C}EK~$^{\ddag}$,\\
 Yuri F. SMIRNOV~$^{\dag\S}$ and Valeriy N.~TOLSTOY~$^{\dag\ddag}$}

\AuthorNameForHeading{R.M.~Asherova, \v{C}. Burd\'{\i}k, M.~Havl\'{\i}\v{c}ek, Yu.F. Smirnov and V.N. Tolstoy}

\Address{$^{\dag}$~Institute of Nuclear Physics, Moscow State University, 119992 Moscow,
Russia}
\EmailD{\href{mailto:raya.acherova@gmail.com}{raya.acherova@gmail.com}, \href{mailto:tolstoy@nucl-th.sinp.msu.ru}{tolstoy@nucl-th.sinp.msu.ru}}

\Address{$^{\ddag}$~Department of Mathematics, Faculty of Nuclear Sciences and Physical
Engineering,\\
\hphantom{$^{\ddag}$}~Czech Technical University in Prague, Trojanova 13, 12000 Prague 2,
Czech Republic}
\EmailD{\href{mailto:burdik@kmlinux.fjfi.cvut.cz}{burdik@kmlinux.fjfi.cvut.cz}, \href{mailto:miloslav.havlicek@fjfi.cvut.cz}{miloslav.havlicek@fjfi.cvut.cz}}

\Address{$^{\S}$~Deceased}

\ArticleDates{Received November 05, 2009, in f\/inal form January 15, 2010;  Published online January 26, 2010}

\Abstract{For the quantum algebra $U_{q}(\mathfrak{gl}(n+1))$ in its reduction on the
subalgeb\-ra $U_{q}(\mathfrak{gl}(n))$ an explicit description of a Mickelsson--Zhelobenko reduction
$Z$-algebra \mbox{$Z_{q}(\mathfrak{gl}(n{+}1),\mathfrak{gl}(n))$} is given in terms of the
generators and their def\/ining relations. \mbox{Using} this $Z$-algebra we describe Hermitian
irreducible representations of a discrete series for the noncompact quantum algebra
$U_{q}(u(n,1))$ which is a real form of $U_{q}(\mathfrak{gl}(n+1))$, namely, an
orthonormal Gelfand--Graev basis is constructed in an explicit form.}

\Keywords{quantum algebra; extremal projector; reduction algebra; Shapovalov form;
noncompact quantum algebra; discrete series of representations; Gelfand--Graev basis}

\Classification{17B37; 81R50}

\renewcommand{\thefootnote}{\arabic{footnote}}
\setcounter{footnote}{0}

\section{Introduction}

In 1950, I.M. Gelfand and M.L. Tsetlin \cite{GT} proposed a formal description of
f\/inite-dimensional irreducible representations (IR) for the compact Lie algebra $u(n)$.
This description is a genera\-li\-zation of the results for $u(2)$ and $u(3)$ to the $u(n)$ case.
It is the following. In the IR space of~$u(n)$ there is a orthonormal basis which is
numerated by the following formal schemes:
\begin{gather*}
\left( \begin{array}{cccccccccc}
m_{1n}&&m_{2n}&&\ldots&&m_{n-1,n}&&m_{nn}\\
&m_{1,n-1}&&m_{2,n-1}&&\ldots&&m_{n-1,n-1}\\
&&\ldots&& \ldots&&\ldots&\\
&&&m_{12}&&m_{22}\\
&&&&m_{11}
\end{array} \right),
\end{gather*}
where all numbers $m_{ij}$ ($1\leq i\leq j\leq n$) are nonnegative integers and
they satisfy the standard inequalities, ``between conditions'':
\begin{gather*}
m_{ij+1}\geq m_{ij}\geq m_{i+1j+1}\qquad{\rm for}\ \ 1\leq i\leq j\leq n-1.
\end{gather*}
The f\/irst line of this scheme is def\/ined by the components of the highest weight of
$u(n)$ IR, the second line is def\/ined by the components of the highest weight of $u(n-1)$
IR and so on.

Later this basis was constructed in many papers (see, e.g.,~\cite{NM,How,AST}) by
using one-step lowering and raising operators.

In 1965, I.M.~Gelfand and M.I.~Graev \cite{GG}, using analytic continuation of the results
for $u(n)$, obtained some results for noncompact Lie algebra $u(n,m)$. It was shown that
some class of Hermitian IR of $u(n,m)$ is characterized by an ``extremal weight''
parametrized by a set of integers $m_{N}=(m_{1N},m_{2N},\ldots, m_{NN})$ $(N=n+m)$ such
that $m_{1N}\geq m_{2N}\geq\cdots \geq m_{NN}$, and by a representation type which is
def\/ined by a partition of $n$ in the sum of two nonnegative integers $\alpha$ and
$\beta$, $n=\alpha+\beta$ (also see~\cite{BR}).

For simplicity we consider the case $u(2,1)$. In this case we have three types of schemes
\begin{gather*}
\left( \begin{array}{ccccc}
& m_{13}&&m_{23}&m_{33}\\
m_{12}&&m_{22}&\\
&m_{11}&
\end{array}\right)\qquad{\rm for}\ \   (\alpha,\beta) = (2,0),
\\ 
\left( \begin{array}{ccccc}
&m_{13}&m_{23}&m_{33}&\\
m_{12}&&&&m_{22}\\
&&m_{11}&
\end{array}\right)\qquad{\rm for}\ \ (\alpha,\beta)=(1,1),
\\
\left(\begin{array}{ccccc}
m_{13}&m_{23}&&m_{33}&\\
&&m_{12}&&m_{22}\\
&&&m_{11}&
\end{array}\right)\qquad{\rm for}\ \ (\alpha,\beta)=(0,2).
\end{gather*}
The numbers $m_{ij}$ of the f\/irst scheme satisfy the following inequalities
\begin{gather*}
m_{12}\geq m_{13}+1,\qquad m_{13}+1\geq m_{22}\geq m_{23}+1,\qquad m_{12}\geq
m_{11}\geq m_{22}.
\end{gather*}
The numbers $m_{ij}$ of the second scheme satisfy the following inequalities
\begin{gather*}
m_{12}\geq m_{13}+1,\qquad m_{33}-1\geq m_{22},\qquad m_{12}\geq m_{11}\geq m_{22}.
\end{gather*}
The numbers of the third scheme satisfy the following inequalities
\begin{gather*}
m_{23}-1\geq m_{12}\geq m_{33}-1,\qquad m_{33}-1\geq m_{22},\qquad m_{12}\geq
m_{11}\geq m_{22}.
\end{gather*}
Construction of the Gelfand--Graev  basis for $u(n,m)$ in terms of one-step lowering and
raising ope\-ra\-tors is more complicated than in the compact case $u(n+m)$.

In 1975, T.J.~Enright and V.S.~Varadarajan \cite{EV} obtained a classif\/ication of discrete
series of noncompact Lie algebras. Later it was proved by A.I.~Molev \cite{Mol} that in
the case of $u(n,m)$ the Gelfand--Graev modules are part of the Enright--Varadarajan modules and
Molev constructed the Gelfand--Graev basis for $u(n,m)$ in terms of the Mickelsson S-algebra~\cite{Mick}.

A goal of this work is to obtain analogous results for the noncompact quantum algebra
$U_q(u(n{,}m))$. Since the general case is very complicated we at f\/irst consider the case
$U_q(u(n{,}1))$. The special case $U_q(u(2,1))$ was considered in~\cite{SK,SKA}. It should
be noted that the principal series representations of $U_q(u(n,1))$ were studied in~\cite{GIK} and a classif\/ication of unitary highest weight modules of $U_q(u(n,1))$ was
considered in~\cite{Gu}.

\section[Quantum algebra $U_{q}(\mathfrak{gl}(N))$ and its noncompact real forms $U_{q}(u(n,m))$ ($n+m=N$)]{Quantum algebra $\boldsymbol{U_{q}(\mathfrak{gl}(N))}$ and its noncompact\\ real forms $\boldsymbol{U_{q}(u(n,m))}$ ($\boldsymbol{n+m=N}$)}
\label{section2}

The quantum algebra $U_{q}(\mathfrak{gl}(N))$ is generated by the Chevalley elements
$q^{\pm e_{ii}}$ $(i=1,\ldots,N)$, $e_{i,i+1}$, $e_{i+1,i}$ $(i=1,2,\ldots,N-1)$
with the def\/ining relations \cite{T1,T2}:
\begin{gather}\label{qa1}
q^{e_{ii}}q^{-e_{ii}} = q^{-e_{ii}}q^{e_{ii}}=1,
\\
\label{qa2}
q^{e_{ii}}q^{e_{jj}}=q^{e_{jj}}q^{e_{ii}},
\\
\label{qa3}
q^{e_{ii}}e_{jk} q^{-e_{ii}}=q^{\delta_{ij}-\delta_{ik}} e_{jk} \qquad(|j-k|=1),
\\
\label{qa4}
[e_{i,i+1} ,e_{j+1,j} ]=\delta_{ij}
\frac{q^{e_{ii}-e_{i+1,i+1}}-q^{e_{i+1,i+1}-e_{ii}}} {q-q^{-1}},
\\
\label{qa5}
[e_{i,i+1} ,e_{j,j+1} ]=0\qquad{\rm for}\ \ |i-j|\geq 2,
\\
\label{qa6}
[e_{i+1,i} ,e_{j+1,j} ]=0\qquad{\rm for}\ \ |i-j|\geq 2,
\\ \label{qa7}
[[e_{i,i+1} ,e_{j,j+1} ]_{q} ,e_{j,j+1} ]_{q} =0\qquad{\rm
for}\ \ |i-j|=1,
\\ \label{qa8}
[[e_{i+1,i} ,e_{j+1,j} ]_{q} ,e_{j+1,j} ]_{q} =0\qquad{\rm
for} \ \ |i-j|=1,
\end{gather}
where $[e_{\beta} ,e_{\gamma} ]_{q} $ denotes the $q$-commutator:
\begin{equation*}
[e_{\beta} ,e_{\gamma} ]_{q} :=e_{\beta} e_{\gamma} -
q^{(\beta,\gamma)}e_{\gamma} e_{\beta}.
\end{equation*}
The def\/inition of a quantum algebra also includes operations of  a comultiplication
$\Delta_{q}$, an antipode $S_{q}$, and a co-unit $\epsilon_{q}$. Explicit
formulas of these operations will not be used in our later calculations and they are not given
here.

Let $\varepsilon_{i}$ ($i=1,2,\ldots,N$) be a dual basis to the Cartan basis $e_{ii}$
($i=1,2,\ldots,N$), $\varepsilon_i (e_{jj} )=(\varepsilon_i,\varepsilon_j)=
\delta_{ij}$. In  terms of the orthonormal basis vectors $\varepsilon_i$ the positive
root system $\Delta_{+}$ of $\mathfrak{gl}(N)$ ($U_{q}(\mathfrak{gl}(N))$) is presented as follows:
\begin{gather*}
\Delta_{+}=\{\varepsilon_i-\varepsilon_j \, |\,1\leq i<j\leq N\},
\end{gather*}
where $\varepsilon_i-\varepsilon_{i+1}$ ($i=1,2,\ldots,N-1$) are the simple roots.

Since for construction of the composite root vectors $e_{ij}:=e_{\varepsilon_i-\varepsilon_j}$
($|i-j|\geq2$) of the quantum algebra $U_{q}(\mathfrak{gl}(N))$ we need to use the
notation of the normal ordering in the positive root system $\Delta_{+}$, we recall this
notation. {\it We say that the system $\Delta_+$ is written in normal $($convex$)$ ordering,
$\vec{\Delta}_+$, if each positive composite root $\varepsilon_i- \varepsilon_j=
(\varepsilon_i-\varepsilon_k)+(\varepsilon_k-\varepsilon_j)$ $(i\leq k\leq j)$ is located
between its components $\varepsilon_i-\varepsilon_k$ and $\varepsilon_k- \varepsilon_j$}.
It means that in the normal ordering system $\vec{\Delta}_+$ we have either
\begin{gather*}
\ldots,\varepsilon_i-\varepsilon_k,\ldots,\varepsilon_i-\varepsilon_j,\ldots,\varepsilon_k-
\varepsilon_j,\ldots,
\end{gather*}
or
\begin{gather*}
\ldots,\varepsilon_k-\varepsilon_j,\ldots,\varepsilon_i-\varepsilon_j,\ldots,\varepsilon_i-
\varepsilon_k,\ldots.
\end{gather*}
There are many normal orderings in the root system $\Delta_{+}=\Delta_{+}
(\mathfrak{gl}(N))$, more than $(N-1)!$ for $N>3$. To be def\/inite, we f\/ix the following
normal ordering (see \cite{T1,T2})
\begin{gather}
\varepsilon_1 -\varepsilon_2 \prec\varepsilon_1 -\varepsilon_3
\prec\varepsilon_2 -\varepsilon_3 \prec\varepsilon_{1} -\varepsilon_{4} \prec
\varepsilon_{2} -\varepsilon_{4} \prec\varepsilon_{3} -\varepsilon_{4}
\prec\cdots\prec \nonumber\\
\varepsilon_{1} -\varepsilon_{k} \prec\varepsilon_{2} -\varepsilon_{k}
\prec\cdots\prec\varepsilon_{k-1} -\varepsilon_{k} \prec\cdots\prec
\varepsilon_{1} -\varepsilon_{N} \prec\varepsilon_{2} -\varepsilon_{N}
\prec\cdots\prec\varepsilon_{N-1} -\varepsilon_{N} .\label{qa13}
\end{gather}
According to this ordering, we determine the composite root vectors $e_{ij} $ for $|i-j|\geq2$ as follows:
\begin{gather}\label{qa14}
e_{ij}  := [e_{ik} , e_{kj} ]_{q^{-1}} ,\qquad
e_{ji}  := [e_{jk} , e_{ki} ]_{q},
\end{gather}
where $1\le i<k<j\le N$. It should be stressed that the structure of the composite root
vectors does not depend on the choice of the index $k$ on the right-hand side of the def\/inition~(\ref{qa14}). In particular, we have{\samepage
\begin{gather}
e_{ij} :=[e_{i,i+1} ,e_{i+1,j} ]_{q^{-1}} =[e_{i,j-1} ,
e_{j-1,j} ]_{q^{-1}} ,\nonumber\\
e_{ji} :=[e_{j,i+1} ,e_{i+1,i} ]_{q} =[e_{j,j-1} ,e_{j-1,i} ]_{q} ,\label{qa15}
\end{gather}
where $2\le i+1<j\le N$.}

Using these explicit constructions and def\/ining relations (\ref{qa1})--(\ref{qa8})
for the Chevalley basis it is not hard to calculate the following relations between the
Cartan--Weyl generators $e_{ij}$ ($i,j=1,2,\ldots, N$):
\begin{gather}\label{qa16}
q^{e_{kk} }e_{ij} q^{-e_{kk} } = q^{\delta_{ki} -\delta_{kj} }
e_{ij}  \qquad(1\le i,j,k\le N),
\\
\label{qa17}
[e_{ij} ,e_{ji} ]=\frac{q^{e_{ii} -e_{jj} }-
q^{e_{jj} -e_{ii} }}{q-q^{-1}}\qquad(1\le i<j\le N),
\\
\label{qa18}
[e_{ij} ,e_{kl} ]_{q^{-1}}=\delta_{jk} e_{il} \qquad(1\le
i<j\le k<l\le N),
\\
\label{qa19}
[e_{ik} ,e_{jl} ]_{q^{-1}} =\big(q-q^{-1}\big)e_{jk} e_{il} \qquad(1\le
i<j<k<l\le N),
\\
\label{qa20}
[e_{jk} ,e_{il} ]_{q^{-1}} =0\qquad(1\le i\le j<k < l\le N),
\\
\label{qa20+}
[e_{ik} ,e_{jk} ]_{q^{-1}} =0\qquad(1\le i < j<k \le N),
\\
\label{qa21}
[e_{kl} ,e_{ji} ]=0\qquad(1\le i<j\le k<l\le N),
\\
\label{qa22}
[e_{il} ,e_{kj} ]=0\qquad(1\le i<j<k<l\le N),
\\\label{qa23}
[e_{ji} ,e_{il} ]=e_{jl} q^{e_{ii} -e_{jj} }\qquad(1\le
i<j<l\le N),
\\
\label{qa24}
[e_{kl} ,e_{li} ]=e_{ki} q^{e_{kk }-e_{ll} }\qquad(1\le
i<k<l\le N),
\\
\label{qa25}
[e_{jl} ,e_{ki} ]=\big(q^{-1}-q\big)e_{kl} e_{ji} q^{e_{jj}-e_{kk}}\qquad(1\le
i<j<k<l\le N).
\end{gather}

If we apply the Cartan involution ($e_{ij}^{\star}=e_{ji}$, $q^{\star}=q^{-1}$) to
the formulas (\ref{qa16})--(\ref{qa25}), we get all relations between the elements of
the Cartan--Weyl basis.

The explicit formula for the extremal projector for $U_q(\mathfrak{gl}(N))$, corresponding
to the f\/ixed normal ordering (\ref{qa13}), has the form \cite{T1,T2}
\begin{gather}
p(U_q(\mathfrak{gl}(N)) = p(U_q(\mathfrak{gl}(N-1))(p_{1N} p_{2N} \cdots
p_{N-2,N} p_{N-1,N} )
\nonumber\\
\phantom{p(U_q(\mathfrak{gl}(N))}{} = p_{12} (p_{13} p_{23} )\cdots(p_{1k} \cdots p_{k-1,k} )\cdots
(p_{1N} \cdots p_{N-1,N} ),\label{qa26}
\end{gather}
where the elements $p_{ij}$ ($1\le i<j\le N$) are given by
\begin{gather}
p_{ij}  = \sum\limits_{r=0}^{\infty}\frac{(-1)^r}{[r]!}
\varphi_{ij,r} e_{ij}^{r}e_{ji}^{r},
\qquad
\varphi_{ij,r}  = q^{-(j-i-1)r}\left\{\prod\limits_{s=1}^{r}
[e_{ii} -e_{jj} +j-i+s]\right\}^{-1}.\label{qa27}
\end{gather}
Here and elsewhere the symbol $[x]$ is given as follows 
\begin{gather*}
[x] = \frac{q^{x}-q^{-x}}{q-q^{-1}}.
\end{gather*}
The extremal projector $p:=p(U_q(\mathfrak{gl}(N))$ satisf\/ies the relations:
\begin{gather}\label{qa29}
e_{i,i+1} p=pe_{i+1,i} =0 \qquad (1\leq i\leq N-1),\qquad p^{2}=p.
\end{gather}
The extremal projector $p$ belongs to the Taylor extension $TU_q(\mathfrak{gl}(N))$ of the
quantum algebras $U_q(\mathfrak{gl}(N))$. The Taylor extension
$TU_q(\mathfrak{gl}(N))$ is an associative algebra generated by formal Taylor series of
the form
\begin{gather*}
\sum_{\{\tilde{r}\},\{r\}}C_{\{\tilde{r}\},\{r\}}\big(q^{e_{11}},\ldots,q^{e_{NN}}\big)
e_{21}^{\tilde{r}_{12}}e_{31}^{\tilde{r}_{13}}e_{32}^{\tilde{r}_{23}}\cdots
e_{N,N-1}^{\tilde{r}_{N-1,N}} e_{12}^{r_{12}}e_{13}^{r_{13}}e_{23}^{r_{23}}\cdots
e_{N-1,N}^{r_{N-1,N}}
\end{gather*}
provided that nonnegative integers $\tilde{r}_{12},\tilde{r}_{13},\tilde{r}_{23},\ldots,
\tilde{r}_{N-1,N}$ and $r_{12},r_{13},r_{23},\ldots,r_{N-1,N}$ are subject to the
constraints
\begin{gather*}
\Bigg|\sum_{i<j}\tilde{r}_{ij}-\sum_{i<j}r_{ij}\Bigg|\leq\mbox{const}
\end{gather*}
for each formal series and the coef\/f\/icients $C_{\{\tilde{r}\},\{r\}}(q^{e_{11}},
\ldots,q^{e_{NN}})$ are rational functions of the $q$-Cartan elements $q^{e_{ii}}$. The
quantum algebra $U_{q}(\mathfrak{gl}(N))$ is a subalgebra  of the
Taylor extension $TU_{q}(\mathfrak{gl}(N))$, $U_{q}(\mathfrak{gl}(N))\subset
TU_{q}(\mathfrak{gl}(N))$.

We consider, on the quantum algebra $U_{q}(\mathfrak{gl}(N))$, two
real forms: compact and noncompact.

The compact quantum algebra $U_{q}(\mathfrak{u}(N))$ can be considered as the quantum
algebra $U_{q}(\mathfrak{gl}(N))$ ($N=n+m$) endowed with the additional Cartan involution
$\star$:
\begin{gather}\label{qa32}
h_{i}^{\star}= h_{i}\qquad{\rm for}\ \ i=1,2,\ldots,N,
\\\label{qa33}
e_{i,i+1}^{\star}= e_{i+1,i},\qquad
e_{i+1,i}^{\star}=e_{i,i+1} \qquad {\rm for}\ \ 1\leq i\leq N-1,
\\\label{qa34}
q^{\star}= q\qquad{\rm or}\qquad q^{\star}=q^{-1}.
\end{gather}
Thus we have two compact real forms: with real $q$ ($q^{\star}=q$) and with circular $q$
($q^{\star}=q^{-1}$). In the case of the circular $q$ the Cartan--Weyl basis $e_{ij}$
($i,j=1,2,\ldots, N$) constructed by formu\-las~(\ref{qa14}) is $\star$-invariant, i.e.\
$e_{ij}^{\star}=e_{ji} $ for all $1\leq i,j\leq N$. In the case of the real $q$ this
Cartan--Weyl basis is not $\star$-invariant, since the basis vectors satisfy the
relations $e_{ij}^{\star}=e_{ji}'$ for $|i-j|\geq2$ where the root vectors $e_{ji}'$ are
obtained from (\ref{qa14}) by the replacement $q^{\pm1}\rightarrow q^{\mp1}$.

It is reasonable to consider the real compact form on the Taylor extension
$TU_{q}(\mathfrak{gl}(N))$. In particular, it should be noted that
\begin{gather}\label{qa35}
p^{\star} = p\qquad{\rm for \ real \ and \ circular} \ \ q.
\end{gather}
This property is a direct consequence of  a uniqueness theorem  for the extremal
projector, which states that {\it equations \eqref{qa29} have a unique nonzero
solution in the space of the Taylor extension $TU_{q}(\mathfrak{gl}(N))$ and this
solution does not depend on the choice of normal ordering and on the replacement
$q^{\pm1}\rightarrow q^{\mp1}$ in formulas \eqref{qa14}}.

The noncompact quantum algebra $U_{q}(\mathfrak{u}(n,m))$ can be considered as the
quantum algebra $U_{q}(\mathfrak{gl}(N))$ ($N=n+m$) endowed with the additional Cartan
involution $*$:
\begin{gather}\label{qa36}
h_{i}^*= h_{i}\qquad{\rm for} \ \ i=1,2,\ldots,N,
\\ \label{qa37}
e_{i,i+1}^* = e_{i+1,i},\qquad e_{i+1,i}^*= e_{i,i+1} \qquad {\rm
for}\ \ 1\leq i\leq N-1,\qquad i\neq n,
\\ \label{qa38}
e_{n,n+1}^* =-e_{n+1,n},\qquad e_{n+1,n}^*=-e_{n,n+1},
\\
\label{qa39}
q^*=q\qquad {\rm or}\qquad q^*=q^{-1}.
\end{gather}
We also have two noncompact real forms: with real $q$ ($q^{\star}=q$) and with circular
$q$ ($q^{\star}=q^{-1}$). Below we will consider the real form
$U_{q}(\mathfrak{u}(n,1))$, i.e.\ the case $N=n+1$.

\section[The reduction algebra $Z_{q}(\mathfrak{gl}(n+1),\mathfrak{gl}(n))$]{The reduction algebra $\boldsymbol{Z_{q}(\mathfrak{gl}(n+1),\mathfrak{gl}(n))}$}\label{section3}

In the linear space $TU_q(\mathfrak{gl}(n+1))$ we separate out a subspace of ``two-sided
highest vectors'' with respect to the subalgebra $U_{q}(\mathfrak{gl}(n))\subset
U_{q}(\mathfrak{gl}(n+1))$, i.e.
\begin{gather*}
\tilde{Z}_{q}(\mathfrak{gl}(n+1),\mathfrak{gl}(n))= \big\{x\in
TU_q(\mathfrak{gl}(n+1))\,\big| \,e_{i,i+1} x=xe_{i+1,i} =0,\ 1\leq i\leq
n -1\big\}.
\end{gather*}
It is evident that if $x\in \tilde{Z}_{q}(\mathfrak{gl}(n+1),\mathfrak{gl}(n))$ then
\begin{gather*}
x= pxp,
\end{gather*}
where $p:=p(U_q(\mathfrak{gl}(n))$. Again, using the annihilation properties of the
projection operator $p$ we have that any vector $x\in
\tilde{Z}_{q}(\mathfrak{gl}(n+1),\mathfrak{gl}(n))$ can be presented in the form of a formal Taylor series
on the following monomials
\begin{gather}\label{z3}
pe_{n+1,1}^{r_{1}'}\cdots e_{n+1,n}^{r_{n}'}e_{n,n+1}^{r_{n}}\cdots e_{1,n+1}^{r_{1}}p.
\end{gather}
It is evident that $\tilde{Z}_{q}(\mathfrak{gl}(n+1),\mathfrak{gl}(n))$ is a subalgebra
in $TU_q(\mathfrak{gl}(n+1))$. We consider a subalgebra $Z_{q}(\mathfrak{gl}(n+1),
\mathfrak{gl}(n))$ generated by f\/inite series on monomials (\ref{z3}).

We set
\begin{gather*}
z_{0} :=p,\qquad z_{i} :=pe_{i,n+1} p,\qquad z_{-i} :=pe_{n+1,i} p
\qquad(i=1,2,\ldots,n).
\end{gather*}

\begin{theorem}
The elements $z_{i}$ $(i=0,\pm1,\pm2,\ldots,\pm n)$ generate the unital associative
algebra $Z_{q}(\mathfrak{gl}(n+1),\mathfrak{gl}(n))$ 
and satisfy the following relations
\begin{gather}\label{z5}
z_{0} z_{i} = z_{i} z_{0}  = z_{i} \qquad{\rm
for}\ \ i=0,\pm1,\pm2,\ldots,\pm n,
\\ \label{z6}
z_{i} z_{-j}  = z_{-j} z_{i}  \qquad{\rm for}\ \ 1\leq i,j\leq
n,\quad i\neq j,
\\
\label{z7}
z_{i} z_{j}  = z_{j} z_{i}  \frac{[\varphi_{ij} +1]}
{[\varphi_{ij} ]} \qquad{\rm for}\ \ 1\leq i<j\leq n,
\\ \label{z8}
z_{-i} z_{-j}  = z_{-j} z_{-i}  \frac{[\varphi_{ij} ]}
{[\varphi_{ij} +1]}\qquad{\rm for}\ \ 1\leq i<j\leq n,
\end{gather}
and
\begin{gather}\label{z9}
z_{i} z_{-i}  = \sum_{j=1}^{n}B_{ij} z_{-j} z_{j} +\gamma_{i} z_{0}
\qquad {\rm for}\ \ i=1,2,\ldots,n,
\end{gather}
where
\begin{gather}\label{z10}
B_{ij}  = -\frac{b_{i}^{-}b_{j}^{+}}{[\varphi_{ij} -1]},
\qquad \gamma_{i}  = [\varphi_{i,n+1} -1]b_{i}^{-},
\\ \label{z11}
b_{i}^{\pm} = \prod_{s=i+1}^{n}\frac{[\varphi_{is} \pm1]}{[\varphi_{is} ]},
\qquad \varphi_{ij}  = e_{ii}-e_{jj}+j-i.
\end{gather}
\end{theorem}

\begin{remark}
The relations (\ref{z5}) state that the element $z_{0} $ is an algebraic unit in
$Z_{q}(\mathfrak{gl}(n+1)$, $\mathfrak{gl}(n))$.
\end{remark}

A proof of the theorem can be obtained by direct calculations using the explicit form of
extremal projector (\ref{qa26}), (\ref{qa27}).

It should be noted that the theorem was proved by V.N.T.\ as early as 1989 but it has not
been published up to now, however, the results of the theorem were used for construction of the
Gelfand--Tsetlin basis for the compact quantum algebra $U_{q}(u(n))$ \cite{T1,T2}.

For construction and study of the discrete  series of the noncompact quantum algebra
$U_{q}(\mathfrak{u}(n{,}1)\!$ we need other relations than~(\ref{z9}). The  system
(\ref{z9}) expresses the elements $z_{i}z_{-i}$ in terms of the elements $z_{-i}z_{i}$
$(i=1,2,\ldots,n)$ but we would like to express the elements $z_{-1}z_{1},
\ldots,z_{-\alpha}z_{\alpha}$, $z_{\alpha+1}z_{-\alpha-1},\ldots,z_{n}z_{-n}$ in terms of
the elements $z_{1}z_{-1},\ldots,z_{\alpha}z_{-\alpha},z_{-\alpha-1}z_{\alpha+1},
\ldots,z_{-n}z_{n}$ for $\alpha=0,1,\ldots,n$.\footnote{In the case $\alpha=0$ we have
relations (\ref{z9}) and for $\alpha=n$ we obtain the system inverse to
(\ref{z9}).}$^{,}\,$\footnote{In Section~\ref{section5} the parameter $\alpha$ will characterize a
representation type of the discrete series.}
These relations are given by the proposition.

\begin{proposition}
The elements $z_{-1}z_{1},\ldots,z_{-\alpha}z_{\alpha},z_{\alpha+1}z_{-\alpha-1},
\ldots,z_{n}z_{-n}$ are expressed in terms of the elements $z_{1}z_{-1},
\ldots,z_{\alpha}z_{-\alpha},z_{-\alpha-1}z_{\alpha+1},\ldots,z_{-n}z_{n}$ by the
formulas
\begin{gather}\label{z12}
z_{-i} z_{i}  = \sum_{j=1}^{\alpha}B_{ij}^{(\alpha)}z_{j} z_{-j} +
\sum_{l=\alpha+1}^{n}B_{il}^{(\alpha)}z_{-l} z_{l} +\gamma_{i}^{(\alpha)}z_{0} \qquad
 (1\leq i\leq\alpha),
\\ \label{z13}
z_{k} z_{-k}  = \sum_{j=1}^{\alpha}B_{kj}^{(\alpha)}z_{j} z_{-j} +
\sum_{l=\alpha+1}^{n}B_{kl}^{(\alpha)}z_{-l} z_{l} +\gamma_{k}^{(\alpha)}z_{0} \qquad
(\alpha+1\leq k\leq n).
\end{gather}
Here
\begin{gather}
B_{ij}^{(\alpha)} =  \frac{b_{i}^{(\alpha)+}b_{j}^{(\alpha)-}}
{[\varphi_{ij} +1]},\qquad  B_{il}^{(\alpha)} = \frac{b_{i}^{(\alpha)+}
b_{l}^{(\alpha)+}}{[\varphi_{il} ]},
\nonumber\\
\gamma_{i}^{(\alpha)} = - [\varphi_{i,n+1} -\alpha]
b_{i}^{(\alpha)+}\qquad {\rm for} \ 1\leq i,j\leq\alpha<l\leq n,
\label{z14}
\\
B_{kj}^{(\alpha)} = - \frac{b_{k}^{(\alpha)-}b_{j}^{(\alpha)-}}
{[\varphi_{kj} ]},\qquad  B_{kl}^{(\alpha)} = -
\frac{b_{k}^{(\alpha)-}b_{l}^{(\alpha)+}}{[\varphi_{kl} -1]},
\nonumber\\
\gamma_{k}^{(\alpha)} = [\varphi_{k,n+1} -\alpha-1]
b_{k}^{(\alpha)-}\qquad{\rm for}\ \ 1\leq j\leq\alpha<k,l\leq n,\label{z15}
\end{gather}
where
\begin{gather}\label{z16}
b_{i}^{(\alpha)\pm} = \left(\prod_{s=1}^{i-1}\frac{[\varphi_{is} \pm1]}
{[\varphi_{is} ]}\right)\left(\prod_{s=\alpha+1}^{n}\frac{[\varphi_{is} ]}
{[\varphi_{is} \pm1]}\right)\qquad (1\leq i\leq\alpha),
\\ \label{z17}
b_{l}^{(\alpha)\pm} = \left(\prod_{s=1}^{\alpha}\frac{[\varphi_{ls} ]}
{[\varphi_{ls} \pm1]}\right)\left(\prod_{s=l+1}^{n}\frac{[\varphi_{ls} \pm1]}
{[\varphi_{ls} ]}\right)\qquad(\alpha+1\leq l\leq n).
\end{gather}
\end{proposition}

\begin{proof}[Scheme of proof] The relations (\ref{z12}) and  (\ref{z13}) with the coef\/f\/icients
(\ref{z14})--(\ref{z17}) can be proved by induction on $\alpha$. For $\alpha=0$ they
coincide with the relations (\ref{z9})--(\ref{z11})\footnote{In this case, the relations
(\ref{z12}) are absent and, moreover, the f\/irst sum in the right side of the relations
(\ref{z13}) is equal to 0 for $\alpha=0$.}. Next we assume that relations
(\ref{z12})--(\ref{z17}) are valid for $\alpha\geq1$ and we extract from (\ref{z13}) the
relation with $k=\alpha+1$ and express in it the term $z_{-\alpha-1}z_{\alpha+1}$ in
terms of the elements $z_{1}z_{-1},\ldots,z_{\alpha+1}z_{-\alpha-1},z_{-\alpha-2}
z_{\alpha+2},\ldots$, $z_{-n}z_{n}$; then this expression is substituted in the right side
of the rest relations (\ref{z12}) and (\ref{z13}) and after some algebraic manipulations
we obtain the relations (\ref{z12})--(\ref{z17}) where $\alpha$ should be replaced by~$\alpha+1$.
\end{proof}

Using (\ref{z6})--(\ref{z8}) and (\ref{z12})--(\ref{z17}) we can prove some power
relations.
\begin{proposition}
The following power relations are valid
\begin{gather}\label{z18}
z_{i}^{r}z_{-j}^{s} = z_{-j}^{s}z_{i}^{r}\qquad {\rm for}\ \ 1\leq i,j\leq
n,\quad i\neq j\quad {\rm and}\quad  r,s\in\mathbb{N},
\\ \label{z19}
z_{i}^{r}z_{j}^{s} = z_{j}^{s}z_{i}^{r}\frac{[\varphi_{ij} +r]!
[\varphi_{ij} -s]!}{[\varphi_{ij} ]![\varphi_{ij} +r-s]!}\qquad {\rm for}\ \  1\leq
i<j\leq n\;\;{\rm and}\quad r,s \in \mathbb{N},
\\ \label{z20}
z_{-i}^{r}z_{-j}^{s} = z_{-j}^{s}z_{-i}^{r}\frac{[\varphi_{ij} ]!
[\varphi_{ij} -r+s]!}{[\varphi_{ij} -r]![\varphi_{ij} +s]!}\qquad {\rm for}\ \
1\leq i<j\leq n\quad {\rm and}\quad r,s \in \mathbb{N},
\\ \label{z21}
z_{-i} z_{i}^{r} = z_{i}^{r-1}\left(\sum_{j=1}^{\alpha}B_{ij}^{(\alpha)}(r)
z_{j} z_{-j} +\!\sum_{l=\alpha+1}^{n}B_{il}^{(\alpha)}(r)z_{-l} z_{l} +
\gamma_{i}^{(\alpha)}(r)z_{0} \right)\qquad  (1\leq i\leq\alpha),\!\!\!
\\
z_{k} z_{-k}^{r} = z_{-k}^{r-1}\left(\sum_{j=1}^{\alpha}B_{kj}^{(\alpha)}(r)
z_{j} z_{-j} +\!\sum_{l=\alpha+1}^{n}B_{kl}^{(\alpha)}(r)z_{-l} z_{l} +
\gamma_{k}^{(\alpha)}(r)z_{0} \right)\nonumber\\
 \qquad \qquad (\alpha\!+\!1\leq k\leq n).\label{z22}
\end{gather}
Here
\begin{gather}
B_{ij}^{(\alpha)}(r) =  \frac{[r]}{[\varphi_{ij} +r]}
b_{i}^{(\alpha)+} (r)  b_{j}^{(\alpha)-},\qquad B_{il}^{(\alpha)}(r) =
\frac{[r]}{[\varphi_{il} +r-1]} b_{i}^{(\alpha)+}(r) b_{l}^{(\alpha)+},
\nonumber\\
\gamma_{i}^{(\alpha)}(r) =  -[r][\varphi_{i,n+1} -\alpha+r-1]
b_{i}^{(\alpha)+}(r)\qquad(1\leq i,j\leq\alpha<l\leq n;\ \ r\in\mathbb{N}),\label{z23}
\\
B_{kj}^{(\alpha)}(r) =  -\frac{[r]}{[\varphi_{kj} -r+1]}
b_{k}^{(\alpha)-} (r) b_{j}^{(\alpha)-},\qquad  B_{kl}^{(\alpha)} =
-\frac{[r]}{[\varphi_{kl} -r]} b_{k}^{(\alpha)-}(r) b_{l}^{(\alpha)+},
\nonumber\\
 \gamma_{k}^{(\alpha)}(r) =  [r][\varphi_{k,n+1} -\alpha-r]
b_{k}^{(\alpha)-}(r)\qquad  (1\leq j\leq\alpha <k,l\leq n;\ \ r\in\mathbb{N}),\label{z24}
\end{gather}
where
\begin{gather}\label{z25}
b_{i}^{(\alpha)+}(r) = \left(\prod_{s=1}^{i-1}\frac{[\varphi_{is} +r]}
{[\varphi_{is} +r-1]}\right)\left(\prod_{s=\alpha+1}^{n}\frac{[\varphi_{is} +r-1]}
{[\varphi_{is} +r]}\right)\qquad (1\leq i\leq\alpha),
\\ \label{z26}
b_{j}^{(\alpha)-} = \left(\prod_{s=1}^{j-1}\frac{[\varphi_{js} -1]}
{[\varphi_{js} ]}\right)\left(\prod_{s=\alpha+1}^{n}\frac{[\varphi_{js} ]}
{[\varphi_{js} -1]}\right)\qquad (1\leq j\leq\alpha),
\\ \label{z27}
b_{k}^{(\alpha)-}(r) = \left(\prod_{s=1}^{\alpha}\frac{[\varphi_{ks} -r+1]}
{[\varphi_{ks} -r]}\right)\left(\prod_{s=k+1}^{n}\frac{[\varphi_{ks} -r]}
{[\varphi_{ks} -r+1]}\right)\qquad(\alpha+1\leq k\leq n),
\\ \label{z28}
b_{l}^{(\alpha)+} = \left(\prod_{s=1}^{\alpha}\frac{[\varphi_{ls} ]}
{[\varphi_{ls} +1]}\right)\left(\prod_{s=l+1}^{n}\frac{[\varphi_{ls} +1]}
{[\varphi_{ls} ]}\right)\qquad (\alpha+1\leq l\leq n).
\end{gather}
\end{proposition}
Here in (\ref{z19}), (\ref{z20}) and thoughtout  in Section~\ref{section4} we use the short notation
of the $q$-factorial $[x+n]!=[x+n][x+n-1]\cdots[x+1][x]!$ instead the $q$-Gamma function,
$[x+n]!\equiv\Gamma_{q}([x+n+1])$.

\begin{proof}[Sketch of proof] The relations (\ref{z18})--(\ref{z20}) are a direct consequence of
the relations (\ref{z6})--(\ref{z8}). Relations (\ref{z21}) and (\ref{z22}) with the
coef\/f\/icients (\ref{z23})--(\ref{z28}) are proved by induction on $r$ using the initial
relations (\ref{z12}) and (\ref{z13}) with the coef\/f\/icients (\ref{z14})--(\ref{z17}) for
$r=1$.
\end{proof}

\section[Shapovalov forms on $Z_{q}(\mathfrak{gl}(n+1),\mathfrak{gl}(n))$]{Shapovalov forms on $\boldsymbol{Z_{q}(\mathfrak{gl}(n+1),\mathfrak{gl}(n))}$}\label{section4}

Let us consider properties of the $Z$-algebra $Z_{q} (\mathfrak{gl}(n+1),
\mathfrak{gl}(n))$ with respect to the involutions~$\star$ (\ref{qa32})--(\ref{qa34}),
and~$*$ (\ref{qa36})--(\ref{qa39}).
\begin{proposition}\label{sf0}
The $Z$-algebra $Z_{q} (\mathfrak{gl}(n+1),\mathfrak{gl}(n))$ is invariant with respect
to the involutions~$\star$ and~$*$, besides
\begin{gather}\label{sf1}
z_{0}^{\star} = z_{0},\qquad z_{i}^{\star}=z_{-i} \qquad {\rm
for}\ \ i=\pm1,\pm2,\ldots,\pm n,
\end{gather}
and
\begin{gather}\label{sf2}
z_{0}^{*} = z_{0},\qquad z_{i}^{*}=-z_{-i} \quad {\rm for}\ \
i=\pm1,\pm2,\ldots,\pm n.
\end{gather}
\end{proposition}

\begin{proof}  Because the extremal projector $p=p(\mathfrak{gl}(n))$ is
$\divideontimes$-invariant, $p^{\divideontimes}=p$ for $\divideontimes=\star,\ast$ (see
the formula (\ref{qa35})), it turns out that
\begin{gather*}
z_{i}^{\divideontimes} = pe_{i,n+1}^{\divideontimes}p,\qquad
z_{-i}^{\divideontimes}=pe_{n+1,i}^{\divideontimes}p\qquad {\rm for}\ \
i=1,2,\ldots,n.
\end{gather*}
If $q$ is circular, then $e_{i,n+1}^{\divideontimes}=\pm e_{n+1,i} $,
$e_{n+1,i}^{\divideontimes}=\pm e_{i,n+1} $, where the plus belongs to the compact case
and the minus belongs to the noncompact case, and we obtain the formulas (\ref{sf1}) and~(\ref{sf2}).

If $q$ is real, then $e_{n+1,i}^{\divideontimes}=\pm e_{i,n+1}'$,
$e_{i,n+1}^{\divideontimes}=\pm e_{n+1,i}'$ where the root vectors $e_{i,n+1}'$ and
$e_{n+1,i}'$ are obtained from (\ref{qa14}) by the replacement $q^{\pm1}\rightarrow
q^{\mp1}$. Let us consider the dif\/ference $z_{-i}^{\divideontimes}-z_{i} =
p(e_{i,n+1}'-e_{i,n+1} )p$ for $1\leq i\leq n$). Substituting here (see the formulas
(\ref{qa15}))
\begin{gather}\label{sf4}
e_{i,n+1}' = e_{in}'e_{n,n+1} -q^{-1}e_{n,n+1} e_{in}',\qquad
e_{i,n+1} =e_{in} e_{n,n+1} -qe_{n,n+1} e_{in},
\end{gather}
and using the annihilation properties of the projector $p$ (see~(\ref{qa29})) we obtain
$z_{-i}^{\divideontimes}-z_{i} =p(e_{in}'-e_{in} )e_{n,n+1} p$. In a similar way,
using explicit formulas of type~(\ref{sf4}) for the generators~$e_{in}'$ and~$e_{in} $, we
obtain $z_{-i}^{\divideontimes}-z_{i} =p(e_{i,n-1}'-e_{i,n-1} )
e_{n-1,n} e_{n,n+1} p$. By proceeding as above, we have
$z_{-i}^{\divideontimes}-z_{i} =p(e_{i,i+1}-e_{i,i+1} )e_{i+1,i+2} e_{i+2,i+3}\cdots
e_{n-1,n} e_{n,n+1} p=0$. In a similar way, we prove that
$z_{i}^{\divideontimes}-z_{-i} =0$.
\end{proof}

The $Z$-algebra $Z_{q} (\mathfrak{gl}(n+1),\mathfrak{gl}(n))$ with the involution
$\star$ is called the compact real form and is denoted by the symbol
$Z_{q}^{(c)}(\mathfrak{gl}(n+1),\mathfrak{gl}(n))$. The noncompact real form on
$Z_{q} (\mathfrak{gl}(n+1),\mathfrak{gl}(n))$ is def\/ined by the involution $*$ and
is denoted by the symbol $Z_{q}^{(nc)}(\mathfrak{gl}(n+1), \mathfrak{gl}(n))$.

Let $p^{(\alpha)}$ be an extremal projector for $Z_{q} (\mathfrak{gl}(n+1),
\mathfrak{gl}(n))$ satisfying the relations
\begin{gather*}
z_{-i} p^{(\alpha)} = p^{(\alpha)}z_{i} \qquad {\rm
for}\ \ i=1,2,\ldots,\alpha,
\\ 
z_{k} p^{(\alpha)} = p^{(\alpha)}z_{-k} \qquad  {\rm
for}\ \ k=\alpha+1,\alpha+2,\ldots,n,
\\ 
[e_{ii}, p^{(\alpha)}] = 0 \qquad {\rm for}\ \ i=1,2,\ldots,n.
\end{gather*}
The extremal projector $p^{(\alpha)}$ depends on the index $\alpha$ that def\/ines what
elements are con\-si\-dered as ``raising'' and what elements are considered as ``lowering'', i.e.
in our case the elements $z_{-1} ,z_{-2} ,\ldots,z_{-\alpha} $,
$z_{\alpha+1} ,\dots, z_{n} $ are raising and the elements
$z_{1} ,z_{2} ,\ldots,z_{\alpha} , z_{-\alpha-1} ,\dots,z_{-n} $ are lowering.
It should be stressed that the ``raising'' and ``lowering'' subsets generate disjoint
subalgebras in $Z_{q} (\mathfrak{gl}(n+1), \mathfrak{gl}(n))$. The operator~$p^{(\alpha)}$ can be constructed in an explicit form.

Let us introduce on $Z_{q}^{(nc)}(\mathfrak{gl}(n+1),\mathfrak{gl}(n))$ the following
sesquilinear Shapovalov form~\cite{Sh}. For any elements $x,y\in
Z_{q}^{(nc)}(\mathfrak{gl}(n+1), \mathfrak{gl}(n))$ we set
\begin{gather}\label{sf8}
B^{(\alpha)}(x,y) = p^{(\alpha)}y^{*}xp^{(\alpha)}.
\end{gather}
Therefore, the Shapovalov form also depends on the index $\alpha$
($\alpha=0,1,2,\ldots,n$). We f\/ix $\alpha$ ($\alpha=0,1,2,\ldots,n$) and for each set of
nonnegative integers $\{r\}=(r_1 ,r_2 ,\ldots,r_n )$ introduce a~vector $v^{(\alpha)}_{\{r\}}$ in the
space $Z_{q}^{(nc)}(\mathfrak{gl}(n+1),\mathfrak{gl}(n))$ by the formula
\begin{gather}\label{sf9}
v^{(\alpha)}_{\{r\}} = z_{\alpha}^{r_\alpha}\cdots
z_{1}^{r_1}z_{-\alpha-1}^{r_\alpha+1}\cdots z_{-n}^{r_n}.
\end{gather}

\begin{theorem}
For each fixed $\alpha$ $(\alpha=0,1,2,\ldots,n)$ the vectors $\{v^{(\alpha)}_{\{r\}}\}$
are pairwise orthogonal with respect to the Shapovalov form \eqref{sf8}
\begin{gather}\label{sf10}
B^{(\alpha)}\big(v^{(\alpha)}_{\{r\}},v^{(\alpha)}_{\{r'\}}\big) = \delta_{\{r\},\{r'\}}
B^{(\alpha)}\big(v^{(\alpha)}_{\{r\}},v^{(\alpha)}_{\{r\}}\big).
\end{gather}
and
\begin{gather}
B^{(\alpha)}\big(v^{(\alpha)}_{\{r\}},v^{(\alpha)}_{\{r\}}\big) =
\Bigg(\prod_{i=1}^{\alpha}\frac{[r_{i}]![\varphi_{i,n+1} -\alpha+r_{i} -1]!}
{[\varphi_{i,n+1} -\alpha-1]!}\prod_{l=\alpha+1}^{n} \frac{[r_{l}]!
[\varphi_{n+1,l} +\alpha+r_{l} ]!} {[\varphi_{n+1,l} +\alpha]!}
\nonumber\\
 \phantom{B^{(\alpha)}\big(v^{(\alpha)}_{\{r\}},v^{(\alpha)}_{\{r\}}\big) = }\times \prod_{1\leq i<j\leq
\alpha}\frac{[\varphi_{ij} +r_{i} -r_{j} ]![\varphi_{ij} -1]!}
{[\varphi_{ij} +r_{i} ]![\varphi_{ij} -r_{j} -1]!}\prod_{\alpha+1\leq k<l\leq
n}\frac{[\varphi_{kl} -r_{k} +r_{l} ]! [\varphi_{kl} -1]!}
{[\varphi_{kl} -r_{k} -1]![\varphi_{kl} +r_{l} ]!}
\nonumber\\
\phantom{B^{(\alpha)}\big(v^{(\alpha)}_{\{r\}},v^{(\alpha)}_{\{r\}}\big) = }
\times \prod_{1\leq i\leq\alpha<l\leq n}
\frac{[\varphi_{il} +r_{i} -1]![\varphi_{il} +r_{l} -1]![\varphi_{il} ]}
{[\varphi_{il} +r_{i} +r_{l} ]![\varphi_{il} -1]!}\biggr)z_{0}^{(\alpha)},\label{sf11}
\end{gather}
where $z_{0}^{(\alpha)}\equiv p^{(\alpha)}$.
\end{theorem}

As a  consequence of this theorem we obtain  that the Shapovalov form is not degenerate on a~subspace of $Z_{q}^{(nc)}(\mathfrak{gl}(n+1),\mathfrak{gl}(n))$, generated by the vectors
of form~(\ref{sf9}).

In the case of  the compact $Z$-algebra $Z_{q}^{(c)}(\mathfrak{gl}(n+1),
\mathfrak{gl}(n))$ the Shapovalov form $B(x,y)$ is def\/ined by formula~(\ref{sf8})
where $\alpha=0$, $p^{(0)}$ is the standard extremal projector of the quantum algebra
$U_{q}(\mathfrak{gl}(n+1))$ and the involution is given by formulas~(\ref{sf1}). It
is not dif\/f\/icult to see that
\begin{gather*}
B(v_{\{r\}},v_{\{r'\}}) = \delta_{\{r\},\{r'\}}B(v_{\{r\}},v_{\{r\}}).
\end{gather*}
where $v_{\{r\}}:= v^{(0)}_{\{r\}}$ and
\begin{gather*}
B(v_{\{r\}},v_{\{r\}}) = (-1)^{\sum\limits_{i=1}^{n}r_{i}}
B^{(0)}\big(v^{(0)}_{\{r\}},v^{(0)}_{\{r\}}\big)
\\
\phantom{B(v_{\{r\}},v_{\{r\}})}{} =  \left(\prod_{l=1}^{n}\frac{[r_{l}]![\varphi_{l,n+1} -1]!}
{[\varphi_{l,n+1} -r_{l} -1]!}\prod_{1\leq k<l\leq
n}\frac{[\varphi_{kl} -r_{k} +r_{l} ]! [\varphi_{kl} -1]!}
{[\varphi_{kl} -r_{k} -1]![\varphi_{kl} +r_{l} ]!}\right)z_{0}^{(0)}.
\end{gather*}

\section[Discrete series of representations for $U_q(u(n,1))$]{Discrete series of representations for $\boldsymbol{U_q(u(n,1))}$}\label{section5}

As in the classical case \cite{Mick} each Hermitian irreducible representation of
the discrete series for the noncompact quantum algebra $U_q(u(n,1))$ is def\/ined uniquely
by some extremal vector~$|xw\rangle$, the vector of extremal weight\footnote{We assume
that the vector $|xw\rangle$ is orthonormal, $\langle xw|xw\rangle=1$.}. This vector
should be the highest vector with respect to the compact subalgebra $U_q(u(n))\oplus
U_q(u(1))$. Since the quantum algebra $U_q(u(1))$ is generated only by one Cartan element
$q^{e_{n+1,n+1} }$, the vector $|xw\rangle$ should be annihilated by the raising
generators $e_{ij} $ ($1\leq i<j\leq n$) of the compact subalgebra $U_q(u(n))$. So the
vector $|xw\rangle$ satisf\/ies the relations
\begin{gather*}
e_{ii} |xw\rangle = \mu_{i}|xw\rangle\qquad(i=1,2,\ldots, n+1),
\\ 
e_{ij} |xw\rangle = 0\qquad (1\leq i<j\leq n),
\end{gather*}
where the weight components $\mu_{i}$ ($i=1,2,\ldots,n)$ are integers subjected to the
condition $\mu_{1}\geq\mu_{2}\geq\cdots\geq\mu_{n}$.
Such weights can be compared with respect to standard lexicographic ordering, namely,
$\mu>\mu'$, where  $\mu=(\mu_{1} ,\mu_{2} ,\ldots,\mu_{n} )$ and
$\mu'=(\mu_{1}^{\prime},\mu_{2}^{\prime},\ldots,\mu_{n}^{\prime})$, if a f\/irst nonvanishing component of
the dif\/ference $\mu-\mu'$ is positive.

The component $\mu_{n+1} $ is also an integer. In the general case of f\/inite-dimensional
irreducible representations of the compact quantum algebra $U_q(u(n))\oplus U_q(u(1))$,
the weights $\mu{=}(\mu_{1} ,\mu_{2}$, $\ldots,\mu_{n} )$ and $\mu_{n+1} $ are not
ordering. If we choose some ordering for these weights, for example,
$(\mu_{1} ,\ldots,\mu_{\alpha} ,\mu_{n+1} ,\mu_{\alpha+1} ,\ldots,\mu_{n} )$,
then such $n+1$-component weights can be compared.

The extremal vector $|xw\rangle$ has the minimal weight $\Lambda^{(\alpha)}_{n+1}:=
(\lambda_{1,n+1} ,\lambda_{2,n+1} ,\ldots,\lambda_{n+1,n+1} )$ where
$\lambda_{i,n+1} :=\mu_{i} $ ($i=1,2,\ldots,\alpha$), $\lambda_{\alpha+1,n+1} :=
\mu_{n+1} $, $\lambda_{l+1,n+1}  :=\mu_{l} $ ($l=\alpha+1,\ldots,n$).
The vector $|\Lambda^{(\alpha)}_{n+1}\rangle:=|xw\rangle$ with the weight
$\Lambda^{(\alpha)}_{n+1}$ satisf\/ies the relations
\begin{gather*}
z_{-i}|\Lambda^{(\alpha)}_{n+1}\rangle = 0,\qquad{\rm for}\ \ i=1,2,\ldots,\alpha,
\\ 
z_{k}|\Lambda^{(\alpha)}_{n+1}\rangle = 0,\qquad{\rm for}\ \ k=\alpha+1,\alpha+2,
\ldots,n.
\end{gather*}
It is evident that any highest weight vector $|\Lambda^{(\alpha)}_{n+1};\Lambda_{n})$
with respect to the compact subalgebra $U_q(u(n))$ has the form
\begin{gather}\label{ds5}
|\Lambda^{(\alpha)}_{n+1};\Lambda_{n}) = z_{\alpha}^{r_\alpha}\cdots
z_{1}^{r_1}z_{-\alpha-1}^{r_{\alpha+1}}\cdots
z_{-n}^{r_n}|\Lambda^{(\alpha)}_{n+1}\rangle.
\end{gather}
Here the integers $\{r\}$ are def\/ined by the weights $\Lambda^{(\alpha)}_{n+1}=
(\lambda_{1,n+1} ,\lambda_{2,n+1} ,\ldots, \lambda_{n+1,n+1} )$, where
$\lambda_{i,n+1} \geq\lambda_{i+1,n+1} $ $(i=1,2,\ldots,n)$, and $\Lambda_{n}=
(\lambda_{1n} ,\lambda_{2n} ,\ldots,\lambda_{nn} )$, where
$\lambda_{in} \geq\lambda_{i+1,n} $ $(i=1,2,\ldots,n-1)$, namely,
\begin{gather*}
r_{i} = \lambda_{in}-\lambda_{i,n+1}\qquad(i=1,\ldots,\alpha),
\\ 
r_{l} = \lambda_{l+1,n+1}-\lambda_{ln}\qquad(l=\alpha+1,\ldots,n).
\end{gather*}
If we would like to calculate the scalar product of the two vectors (\ref{ds5}) then using
the results for the Shapovalov form (\ref{sf10}), (\ref{sf11}) we obtain
\begin{gather*}
\big(\Lambda_{n};\Lambda^{(\alpha)}_{n+1}|\Lambda^{(\alpha)}_{n+1};\Lambda_{n}'\big) =
\delta_{\Lambda_{n},\Lambda_{n}'}\big(\Lambda_{n};\Lambda^{(\alpha)}_{n+1}|
\Lambda^{(\alpha)}_{n+1};\Lambda_{n}\big),
\\ 
\big(\Lambda_{n};\Lambda^{(\alpha)}_{n+1}|\Lambda^{(\alpha)}_{n+1};\Lambda_{n}\big) =
B^{(\alpha)}\big(v^{(\alpha)}_{\{r\}},v^{(\alpha)}_{\{r\}}\big)\Bigr|_{\Lambda^{(\alpha)}_{n+1}},
\end{gather*}
where the symbol $|_{\Lambda^{(\alpha)}_{n+1}}$ means that we specialize the Shapovalov form
(\ref{sf11}) for the extremal weight~$\Lambda_{n+1}$, that is we replace the Cartan
elements $e_{ii}$, $e_{jj} $ in the functions~$\varphi_{ij}$ by the corresponding compo\-nents~$\lambda_{i,n+1}$, $\lambda_{j,n+1}$ and~$z_0$ by~$1$.

From the condition that
\begin{gather*}
\big(\Lambda_{n};\Lambda^{(\alpha)}_{n+1}|\Lambda^{(\alpha)}_{n+1};\Lambda_{n}\big) > 0
\end{gather*}
we f\/ind all admissible highest weights $\Lambda_{n}$ of the compact subalgebra
$U_q(u(n))$. The result is formulated as the theorem 
\begin{theorem}
$1)$  Every Hermitian irreducible representation of the discrete series for the noncompact
quantum algebra $U_q(u(n,1))$ with the extremal weight $\Lambda^{( \alpha)}_{n+1}=
(\lambda_{1,n+1} ,\ldots, \lambda_{n+1,n+1} )$, where the integers $\lambda_{i,n+1}$
satisfy the inequalities $\lambda_{i,n+1} \geq\lambda_{i+1,n+1} $ $(i=1,2,\ldots,n)$,
under the restriction $U_q(u(n,1))\downarrow U_q(u(n))$ contains all multiplicity free
irreducible representations of the compact subalgebra $U_q(u(n))$ with the highest weights
$\Lambda_{n}=(\lambda_{1n} ,\lambda_{2n} ,\ldots, \lambda_{nn} )$ satisfying the
conditions:
\begin{gather}
 \lambda_{1n}\geq\lambda_{1,n+1}\geq\lambda_{2,n}\geq\lambda_{2,n+1}\geq\cdots\geq
\lambda_{\alpha n}\geq\lambda_{\alpha,n+1},
\nonumber\\
 \lambda_{\alpha+2,n+1}\geq\lambda_{\alpha+1,n}\geq\lambda_{\alpha+3,n+1}\geq\cdots\geq\lambda_{n+1,n+1}\geq
\lambda_{nn}.\label{ds11}
\end{gather}

$2)$ The vectors
\begin{gather*}
|\Lambda^{(\alpha)}_{n+1};\Lambda_{n} \rangle = F^{(\alpha)}_{-}(\Lambda_{n};
\Lambda^{(\alpha)}_{n+1})\,|\, \Lambda^{(\alpha)}_{n+1}\rangle,
\end{gather*}
where the ``lowering'' operator $F^{(\alpha)}_{-}(\Lambda_{n} ;
\Lambda^{(\alpha)}_{n+1})$ is given by
\begin{gather}
F^{(\alpha)}_{-}\big(\Lambda_{n} ;\Lambda^{(\alpha)}_{n+1}\big) =
N^{(\alpha)}\big(\Lambda_{n} ;\Lambda^{(\alpha)}_{n+1}\big)  z_{\alpha}^{\lambda_{\alpha
n}-\lambda_{\alpha,n+1}}\cdots z_{1}^{\lambda_{1n}-\lambda_{1,n+1}}
\nonumber\\
 \phantom{F^{(\alpha)}_{-}\big(\Lambda_{n} ;\Lambda^{(\alpha)}_{n+1}\big) =}{}
 \times z_{-\alpha-1}^{\lambda_{\alpha+2,n+1}-\lambda_{\alpha+1,n}}\cdots
z_{-n}^{\lambda_{n+1,n+1}-\lambda_{nn}},\label{ds13}
\end{gather}
for all highest wights $\Lambda_{n}=(\lambda_{1n} ,\lambda_{2n} ,\ldots,
\lambda_{nn} )$ constrained by the conditions \eqref{ds11} form an orthonormal basis in
the space of the highest vectors with respect to the compact subalgebra $U_q(u(n))$. Here
in \eqref{ds13} the normalized factor $N^{(\alpha)}(\Lambda_{n} ;
\Lambda^{(\alpha)}_{n+1})$ is given as follows:
\begin{gather*}
N^{(\alpha)}\big(\Lambda_{n} ;\Lambda^{(\alpha)}_{n+1}\big) = \big(\Lambda_{n};
\Lambda^{(\alpha)}_{n+1}|\Lambda^{(\alpha)}_{n+1};\Lambda_{n}\big)^{-\frac{1}{2}}
 \\
 \phantom{N^{(\alpha)}\big(\Lambda_{n} ;\Lambda^{(\alpha)}_{n+1}\big)}{} =
\Bigg\{\prod_{i=1}^{\alpha}\frac{[l_{i,n+1} -l_{\alpha+1,n+1} -2\alpha+n-1]!}
{[l_{in} -l_{i,n+1} ]![l_{i,n} -l_{\alpha+1,n+1} -2\alpha+n-1]!}
\\
 \phantom{N^{(\alpha)}\big(\Lambda_{n} ;\Lambda^{(\alpha)}_{n+1}\big)}{}
 \times \prod_{l=\alpha+1}^{n}
\frac{[l_{\alpha+1,n+1} -l_{l+1,n+1} +2\alpha-n-1]!}
{[l_{l+1,n+1} -l_{ln} -1]![l_{\alpha+1,n+1} -l_{ln}+2\alpha-n]!}\bigg.
\\
\phantom{N^{(\alpha)}\big(\Lambda_{n} ;\Lambda^{(\alpha)}_{n+1}\big)}{}
 \times \prod_{1\leq i<j\leq
\alpha}\frac{[l_{in} -l_{j,n+1} ]![l_{i,n+1} -l_{jn} -1]!}
{[l_{in} -l_{jn} ]![l_{i,n+1} -l_{j,n+1}-1]!}
\\
\phantom{N^{(\alpha)}\big(\Lambda_{n} ;\Lambda^{(\alpha)}_{n+1}\big)}{}
 \times \prod_{\alpha+1\leq k<l\leq n}\frac{[l_{kn} -l_{l+1,n+1} -2]!
[l_{k+1,n+1} -l_{ln} +1]!}{[l_{kn} -l_{ln} ]! [l_{k+1,n+1} -l_{l+1,n+1} -1]!}
\\
\phantom{N^{(\alpha)}\big(\Lambda_{n} ;\Lambda^{(\alpha)}_{n+1}\big)}{}  \times  \prod_{1\leq i\leq\alpha<l\leq n}
\frac{[l_{in} -l_{ln}]![l_{i,n+1} -l_{l+1,n+1}-2]!}
{[l_{in} -l_{l+1,n+1}-2]![l_{i,n+1} -l_{ln}-1]![l_{i,n+1} -l_{l+1,n+1}-1]}
\Bigg\}^{\frac{1}{2}}
\end{gather*}
$(l_{sr}:=\lambda_{sr}-s$ for $s=1,2,\ldots,r$; $r=n,n+1)$.
\end{theorem}

The f\/irst part of the theorem coincides with the classical Gelfand--Graev case
\cite{GG,BR} for the noncompact Lie algebra $u(n,1)$. Using analogous construction of the
Gelfand--Tsetlin basis for the compact quantum algebra $U_q(u(n))$ \cite{T1} we obtain a
$q$-analog of the Gelfand--Graev--Tsetlin basis for $U_q(u(n,1))$. Namely, in the $
U_q(u(n,1))$-module with the extremal weight ${\Lambda^{(\alpha)}_{n+1}}$ there is an
orthogonal Gelfand--Graev--Tsetlin basis consisting of all vectors of the form
\begin{gather*}
\bigl|\Lambda\rangle  := \left|
\begin{array}{c}
\Lambda^{(\alpha)}_{n+1}\\
\Lambda_{n}\\
\dots\\
\Lambda_{2}\\
\Lambda_{1}\\
\end{array}\right\}
= F_{-}(\Lambda_1;\Lambda_2)F_{-}(\Lambda_2;\Lambda_3)\cdots
F_{-}(\Lambda_{n-1};\Lambda_n)|\Lambda^{(\alpha)}_{n+1};\Lambda_{n} \rangle,
\end{gather*}
where $\Lambda_j=(\lambda_{1j},\lambda_{2j},\ldots,\lambda_{jj})$ $(j=1,2,\dots,n)$ and
the numbers $\lambda_{ij}$ satisfy 
the standard ``between conditions'' for the quantum algebra $U_q(u(n))$, i.e.
\begin{gather*}
\lambda_{i,j+1}\geq\lambda_{ij}\geq\lambda_{i+1,j+1}\qquad{\rm for}\ \ 1\leq i\leq j\leq
n-1.
\end{gather*}
The lowering operators $F_{-}(\Lambda_k;\Lambda_{k+1})$ ($k=1,2,\ldots,n-1$) are given by
(see \cite{T1, T2})
\begin{gather*}
F_{-}(\Lambda_k;\Lambda_{k+1}) = N(\Lambda_k;\Lambda_{k+1})
p(U_q(u(k)))\prod_{i=1}^{k}(e_{k+1i})^{\lambda_{ik+1}-\lambda_{ik}} ,
\\
N(\Lambda_k;\Lambda_{k+1}) =  \Bigg\{\prod\limits_{i=1}^{k}
\frac{[l_{ik}-l_{k+1,k+1}-1]!}{[l_{i,k+1}-l_{ik}]! [l_{i,k+1}-l_{k+1,k+1}-1]!}
\\
\phantom{N(\Lambda_k;\Lambda_{k+1}) =}{} \times \prod\limits_{1\le i<j\le k}\frac{[l_{i,k+1}-l_{jk}]!
[l_{ik}-l_{j,k+1}-1]!}{[l_{ik}-l_{jk}]![l_{i,k+1}-l_{j,k+1}-1]!}\bigg\}^{\frac{1}{2}},
\end{gather*}
where $l_{ij}:=\lambda_{ij}-i$ for $1\leq i\leq j\leq n-1$. This explicit construction
allows one to obtain formulas for  actions of the $U_q(u(n,1))$-generators. These results will
be presented elsewhere.

\section{Summary}\label{section6}

Thus, we obtain the explicit description of the Hermitian irreducible representations of the discrete series for
the noncompact quantum algebra $U_q(u(n,1))$ by the reduction $Z$-algebras for
description of which we used the standard extremal projectors.

Next step: to obtain analogous results for $U_q(u(n,2))$. For this aim we
need to construct the extremal projector $p^{(\alpha)}$ which is expressed in terms of the
$Z$-algebra $Z_{q} (\mathfrak{gl}(n+1),\mathfrak{gl}(n))$.

Final aim: to consider the general case $U_q(u(n,m))$. In this case, extremal projectors
of new type will be used.

\subsection*{Acknowledgments}
The paper has been supported by grant RFBR-08-01-00392 (R.M.A., V.N.T.) and by grant
RFBR-09-01-93106-NCNIL-a (V.N.T.). The f\/ifth author would like to thank Department of
Mathematics, Faculty of Nuclear Sciences and Physical Engineering, Czech Technical
University in Prague for hospitality.

\pdfbookmark[1]{References}{ref}
\LastPageEnding

\end{document}